\documentclass[11pt]{article}

\usepackage{amsmath,amsfonts,amssymb,amsthm,bbm}
\usepackage{graphics,epsfig}
\usepackage{color}
\usepackage{array}
\usepackage{ulem}
\usepackage{mathrsfs}
\usepackage{extarrows}
\usepackage{tikz}
\usepackage{indentfirst}

\newtheorem{theorem}{Theorem}[section]
\newtheorem{proposition}{Proposition}[section]

\newtheorem{remark}{\noindent\mbox{Remark}}[section]

\numberwithin{equation}{section}

\allowdisplaybreaks

\renewcommand\today{\number\year/\number\month/\number\day}

\def\qed{\hfill$\Box$\medskip}

\def\z{{\mathbb Z}}
\def\p{{\mathbb P}}
\def\r{{\mathbb R}}
\def\n{{\mathbb N}}
\def\e{{\mathbb E}}

\begin{document}

\noindent\makebox[60mm][1]{\tt {\large Version:~\today}}

\bigskip
%%%%%%%%%%%%%%%%%%%%%%%%%%%%%%% Title %%%%%%%%%%%%%%%%%%%%%%%%%%%%%%%%%%%%%%\\
\noindent{
{\Large\bf  Scaling limit theorems for the $\kappa$-transient random walk in random and non-random environment
%Discrete approximation for Diffusion process in a Brownian environment with drift
}\footnote{
\noindent  The project is partially supported by the National
Natural Science Foundation of China (Grant No. 11131003).
}
\\

\noindent{%\normalsize\sf
Wenming Hong\footnote{ School of Mathematical Sciences
\& Laboratory of Mathematics and Complex Systems, Beijing Normal
University, Beijing 100875, P.R. China. Email: wmhong@bnu.edu.cn} ~  Hui Yang\footnote{ School of Mathematical Sciences
\& Laboratory of Mathematics and Complex Systems, Beijing Normal
University, Beijing 100875, P.R. China. Email: yanghui2011@mail.bnu.edu.cn} ~

\noindent{%\normalsize\sf
(Beijing Normal University)
}

}

\vspace{0.1 true cm}

%%%%%%%%%%%%%%%%%%%%%%%%%%%%%% Abstrct %%%%%%%%%%%%%%%%%%%%%%%%%%%%%%%%%%%%
\begin{center}
\begin{minipage}[c]{12cm}
\begin{center}\textbf{Abstract}\end{center}
\bigskip
 Kesten et al.( 1975) proved the stable law for the transient RWRE (here we refer it as the
$\kappa$-transient RWRE). After that,   some similar interesting properties have also been revealed for its continuous counterpart, the diffusion proces in a Brownian environment with drift $\kappa$.
In the present paper we will investigate the connections between these two kind of models, i.e., we will construct  a sequence of the $\kappa$-transient RWREs and  prove it convergence to the diffusion proces in a Brownian environment with drift $\kappa$ by proper scaling. To this end, we need a counterpart convergence for the $\kappa$-transient random walk in non-random environment, which is interesting itself.

\mbox{}\textbf{Keywords:}\quad Random walk, Random environment, Diffusion process, Brownian motion with drift. \\
\mbox{}\textbf{Mathematics Subject Classification}:  Primary 60F17; secondary 60G50; 60J60;

\end{minipage}
\end{center}

\bigskip
\bigskip

%%%%%%%%%%%%%%%%%%%%%%%%%%%Section 1%%%%%%%%%%%%%%%%%%%%%%%%%%%%%%%%%%%%%%%%%%
\section{ Introduction }
       \label{s1:intro}

%%%%%%%%%%%%%%%%%%%%%%%%%%%%Random walk in non-Random Environment%%%%%%%%%%%%%%%%%%%
%%%%%%%%%%%%%%%%%%%%%%%%%%%%Random walk in Random Environment%%%%%%%%%%%%%%%%%%%%%%%%%

% \subsection{ $\kappa$-Transient Random Walk In Random Environment }           \label{s1.2:random}

 Let $\omega := \{\omega_{i}\}_{i\in\z}$ be a sequence of independent and identically distributed random variables taking values in $(0,1)$, which is labeled as ``random environment". For any realization of $\omega$, we  define a Markov chain $\{Z_{n}\}_{n\geq0}$ with $Z_{0} = 0$ and for any $n \geq 0$ and $i \in \z$,
\begin{equation} \label{rwre}
\p \left( Z_{n+1} = j ~|~Z_{n} = i; \omega \right) = \left\{\begin{array}{lll} \omega_{i}, & {\rm if}
\ j = i+1,
\\
 1-\omega_{i}, &  {\rm if} \ j = i-1,
\\
0, &{\rm otherwise},
\end{array}
\right.
\end{equation}
The process $\{Z_{n}\}_{n\geq0}$ is a so-called random walk in random environment (RWRE, for simplicity), we refer the readers to  Zeitouni (\cite{ze04}, 2004) for more details. We denote $P$   as the law  of the environment, $P_{\omega}$ and  $\p$  as the quenched and annealed probability   of the RWRE respectively, and $E$,   $E_{\omega}$ and $\e$ as the corresponding  expectation.

In the famous paper, Kesten et al.(\cite{K75}, 1975) proved the stable law for the transient RWRE. The main assumptions on the random environment $\omega$  is, (denote $\rho_{i}=\frac{1-\omega_{i}}{\omega_{i}}$, $i \in \z$)
\begin{equation}\label{rwrec}
\begin{cases}
\{\omega_{i}\}_{i\in\z} ~is~ a~ sequence~ of~ i.i.d.~ random~ variables~ taking~ values~ in~ (0,1),\\
E(\log \rho_{0})<0, \\
\exists ~\kappa>0,~s.t.~E\rho_{0}^{\kappa}=1.
\end{cases}
\end{equation}
The condition $E(\log \rho_{0})<0$ indicates the RWRE is transient, and intuitively the $\kappa$ in $E\rho_{0}^{\kappa}=1$ is a ``measurement" for the transient in some sense, and we call it  the $\kappa$-transient RWRE. In addition of the stable law (see \cite{K75}, 1975),  many other interesting results on this $\kappa$-transient RWRE has been obtained, for example, the large deviation principles have been proved by
Greven  and den Hollander (\cite{G94}, 1994), Dembo et al (\cite{DPZ96}, 1996) and Comets et al (\cite{C00}, 2000); see also
 Enriquez et al (\cite{ESZ09}, 2009); etc..

\indent
On the other hand, as a continuous version of the  RWRE, the diffusion process with random potential
has been considered in the literature.
%( Schumacher \cite{S85}).%and Brox \cite{B86}.
Let $V = (V(x), x \in \r)$ be a stochastic process defined on $\r$ with $V(0)=0$. Define $X_{V}=(X_{V}(t), t\geq0)$ a diffusion process in random potential $V$ as an informal solution of the stochastic differential equation
\begin{equation}\label{CX}
\begin{cases}
dX_{V}(t) = dB(t) - \frac{1}{2} V^{\prime}( X_{V}(t) )dt, \\
X_{V}(0) = 0.
\end{cases}
\end{equation}
where $(B(t), t\geq0)$ is a standard Brownian motion independent of $V$.
Strictly speaking, instead of writing the formal derivation of $V$, the process $X_{V}$ can be considered as a diffusion process with generator
\begin{equation*}
\frac{1}{2} e^{V(x)} \frac{d}{dx} \left( e^{-V(x)} \frac{d}{dx} \right).
\end{equation*}
It is called {\it diffusion proces in a Brownian environment with drift $\kappa\geq 0$} if
\begin{equation}\label{PV}
 {V}(x)=\sigma \cdot W(x)-\frac{\kappa}{2}x, ~~~x \in \r,
\end{equation}
where $\{W(x); x\in \r\}$ is a two-sided Brownian motion with $W(0)=0$.

(1) When $\kappa=0$, $X_{V}$ is also called {\it Brox} diffusion, which was discussed by (Schumacher \cite{S85} and Brox \cite{B86}). They showed that $X_{V}$ exhibits the same asymptotic behavior as Sinai's random walk (i.e., $E(\log \rho_{0})=0$), and the {\it Brox} diffusion can be viewed as a continuous version of the Sinai's random walk.  Indeed, Seignourel (\cite{S00}, 2000) proved that   a sequence of Sinai's random walk convergence to the  Brox's diffusion by scaling.

(2)When $\kappa>0$, $X_{V}$ exhibits the similar asymptotic behavior as the $\kappa$-transient RWRE (Kesten et al.\cite{K75}, 1975) such as the stable law (Kawazu and Tanaka \cite{K97},  and Tanaka \cite{T95}) and the large deviations (Hu et al. \cite{H99} and  Taleb \cite{T01}), etc.,  and the  diffusion $X_{V}$ can be viewed as a continuous version of the $\kappa$-transient RWRE in some sense. It is naturally to ask the corresponding scaling limit relations, similar as for the situation $\kappa=0$ by Seignourel (\cite{S00}, 2000).

\

The purpose of the present paper is to construct  a sequence of the $\kappa$-transient RWREs and to prove it convergence to the diffusion proces in a Brownian environment with drift $\kappa$ by proper scaling. To this end, we need a counterpart convergence for the random walk in non-random environment, which is interesting itself.

\section{ \it $\beta$-transient random walk in non-random environment }
          \label{s1.1:non-random}

$\{X_{i}, i\in \n^{+}\}$ is a sequence of i.i.d. random variables taking values in $\r$ whose distribution denoted by $\mu$, and satisfies the following assumptions:
\begin{equation}\label{rwc}
\begin{cases}
1.~EX_{1}<0,\\
2.~\exists ~ \beta > 0 ~s.t.~ Ee^{\beta X_{1}}=1,\\
3.~\mu ~has ~compact ~support.
\end{cases}
\end{equation}
Define a random walk  $\{S_{n}\}_{n\geq0}$ with $S_{0}=0$ and for $n\geq1$,
\begin{equation*}
S_{n}=\sum_{i=1}^{n}X_{i}.
\end{equation*}
Conditions in (\ref{rwc}) states that the random work is transient to $-\infty$ with non-zero speed, and the value of $\beta$   ``measure" the transient  of the random walk in some sense. It is worthwhile to point out that,  by the Donsker's invariant principle with the usual schedule, one can not obtain the limit process as the Brownian motion with drift. To this end, we must construct a sequence of random walks with ``measure change".  Borisov et al (\cite{BN12}) considered the transient simplest random walk case. Here we consider the general situation under the condition (\ref{rwc}),  with the ``transient measurement" $\beta$. A key step is to find the suitable ``measure change" as the following,

\

$\clubsuit$ {\it step 1.   Define $\mu^{(m)}$ }

Starting from $\mu$, which satisfies the condition (\ref{rwc}). For each $m\geq1$, $\mu^{(m)}$, is defined as following, for any continuous bounded function $\varphi$,
\begin{flalign}\label{distri}
\int_{\r} \varphi(x) \mu^{(m)}(dx)  =& c^{-1}\int_{0}^{\infty}\mu(dv)\int^{0}_{-\infty}\mu(du)(e^{\beta v}-e^{\beta u}) \cdot \nonumber \\
&\left[ \left(\frac{-v}{u-v} - \frac{\beta}{2\sqrt{m}(u-v)} \right) \varphi(u) + \left(\frac{u}{u-v} + \frac{\beta}{2\sqrt{m}(u-v)} \right) \varphi(v) \right],
\end{flalign}
where $c:= \int_{0}^{\infty}(e^{\beta x}-1)\mu(dx)=\int^{0}_{-\infty}-(e^{\beta x}-1)\mu(dx) > 0$. Write
\begin{equation}\label{sigma}
\sigma^{2} = -c^{-1}\int_{0}^{\infty}\mu(dv)\int^{0}_{-\infty}\mu(du)(e^{\beta v}-e^{\beta u})uv.
\end{equation}
By condition 3 in (\ref{rwc}), we know $0<\sigma^{2}<\infty$. It is easy to check that $\mu^{(m)}$ is a probability measure for each $m\geq 1$.

$\clubsuit$ {\it step 2.   Define $\{X_{i}^{(m)}, i\in \n^{+}\}$ }

For each $m\geq1$, define a sequence of i.i.d. random variables $\{X_{i}^{(m)}, i\in \n^{+}\},  $ whose distribution is given by $\mu^{(m)}$, which is defined in (\ref{distri}). The corresponding random walk $\{ S_{n}^{(m)}, n \geq 0 \}$ is defined as  $S_{0}^{(m)}=0$, and
\begin{equation*}
S_{n}^{(m)}=\sum_{i=1}^{n}X_{i}^{(m)}.
\end{equation*}
\begin{proposition}\label{th1}{\it Under the assumption (\ref{rwc}), we have
as $m\rightarrow \infty$,
\begin{equation}\label{bmd}
\left\{  \frac{1}{\sqrt{m}}S_{[mt]}^{(m)},~ t \in \r^{+} \right\} \xlongrightarrow[~~~~~]{} \left \{ H(t)=\sigma \cdot B(t) - \frac{\beta}{2}t,~ t\in \r^{+} \right \}
\end{equation}
in distribution in $\mathcal{D}[0,+\infty)$, Where $\sigma$ is defined in (\ref{sigma}) and  $(B(t), t\geq0)$ is a standard Brownian motion.
}
\end{proposition}
\begin{remark} If $\{S_{n}\}_{n\geq0}$ is a simple random walk and satisfies the conditions (\ref{rwc}), then $\mu$ the distribution of $X_{i}$ is,
\begin{flalign*}
&\mu(1)=\frac{1}{e^{\beta}+1}, \\
&\mu(-1)=\frac{e^{\beta}}{e^{\beta}+1}.
\end{flalign*}
It is easy to calculate from (\ref{distri}),
\begin{flalign*}
&\mu^{(m)}(1)=\frac{1}{2}+\frac{\beta}{4\sqrt{m}},\\
&\mu^{(m)}(-1)=\frac{1}{2}-\frac{\beta}{4\sqrt{m}},
\end{flalign*}
which is exactly the random walks constructed in the proof of Theorem 2 in \cite{BN12}. For the general situation, the construction of the $\mu^{(m)}$ in (\ref{distri}) is stimulated by the Skorokhod's representation theorem (see for example page 399 in \cite{Dur}).
\end{remark}

\noindent
{\it Proof of Proposition \ref{th1}.}
 Write
\begin{equation}\label{RW}
H^{(m)}(t)=
\begin{cases}
\frac{1}{\sqrt{m}}S_{k}^{(m)},~ \mbox{if}~ t=\frac{k}{m}, ~k\in \n, \\
\mbox{linear on} ~[\frac{k}{m},\frac{k+1}{m}],~ \mbox{for} ~k\in \n .                                                          \end{cases}
\end{equation}
it is enough to prove the weak convergence of $H^{(m)}(\cdot)$ to $H(\cdot)$ in $\mathcal{C}[0,\infty)$. To this end, we just need to illustrate the finite dimensional distributions of $ H^{(m)}(\cdot)$ converge to the corresponding finite-dimensional distributions of $H(\cdot)$ as $m \rightarrow \infty$ and the relatively compact of $\{H^{(m)}(\cdot)\}_{m\geq1}$ (  Theorem 7.8 of Chapter 3 in (\cite{Kurtz})).

\

(1){\it convergence of the finite dimensional distributions of $ H^{(m)}(\cdot)$}.

By Theorem 6.5 of Chapter 1 and Corollary 8.5 of Chapter 4 in (\cite{Kurtz}), we just need to verify that,
\begin{equation}\label{generator}
mA_{m}f \rightarrow Lf
\end{equation}
uniformly on compact subsets of $\r$, where $L=\frac{1}{2}\sigma^{2}\frac{d^{2}}{(dx)^{2}} + (-\frac{\beta}{2})\frac{d}{dx}$ is the generator of $H(\cdot)$ and
\begin{equation*}
A_{m}f(x)=\int(f(y)-f(x))\Pi_{m}(x,dy),
\end{equation*}
for all $f\in C_{0}^{\infty}(\r)$. $\Pi_{m}(x,\Gamma)$ is transition probability of  function $H^{(m)}$, defined as following,
\begin{flalign}\label{tp}
\Pi_{m}(x,\Gamma) &:= P \left( H^{(m)}\left(\frac{k+1}{m}\right) \in  \Gamma~|~H^{(m)}\left(\frac{k}{m}\right) = x \right)\nonumber\\
&= P \left( \frac{1}{\sqrt{m}}S_{k+1}^{(m)} \in  \Gamma~|~ \frac{1}{\sqrt{m}}S_{k}^{(m)} = x \right) \nonumber \\
&=P \left( \left( x + \frac{1}{\sqrt{m}} X_{k+1}^{(m)}\right) \in \Gamma \right) =\int_{x+\frac{1}{\sqrt{m}}z \in \Gamma}\mu^{(m)}(dz).
\end{flalign}
To prove (\ref{generator}), it is enough to check the following three conditions (\ref{c1}), (\ref{c2}) and (\ref{c3})
(  by {  Lemma 11.2.1} (\cite{V09}))
\begin{flalign}
&\lim_{m\rightarrow\infty}\sup_{|x|\leq R}|a_{m}(x)-\sigma^{2}|=0,\label{c1}\\
&\lim_{m\rightarrow\infty}\sup_{|x|\leq R}|b_{m}(x)-(-\frac{\beta}{2})|=0,\label{c2}\\
&\lim_{m\rightarrow\infty}\sup_{|x|\leq R}\Delta_{m}^{\varepsilon}=0, ~~\varepsilon>0\label{c3}.
\end{flalign}
for all $R>0$, where
\begin{flalign*}
&a_{m}(x)=m\int_{|y-x|\leq1}(y-x)^{2}\Pi_{m}(x,dy),\\
&b_{m}(x)=m\int_{|y-x|\leq1}(y-x)\Pi_{m}(x,dy),\\
&\Delta_{m}^{\varepsilon}=m \Pi_{m}(x,\r \backslash B(x,\varepsilon)),~~\varepsilon>0.
\end{flalign*}
Here, for large enough $m$, by the compact support of the distribution $\mu$ and the expression of $\Pi_{m}(x,\Gamma)$ in (\ref{tp}), it is easy to  calculate that
\begin{flalign*}
a_{m}(x)
&=m\int_{|y-x|\leq1}(y-x)^{2}\Pi_{m}(x,dy),\\
&=m\int_{| \frac{1}{\sqrt{m}}z | \leq 1}\left( \frac{1}{\sqrt{m}}z \right)^{2}\mu^{(m)}(dz),\\
&=\int_{\r}z^{2}\mu^{(m)}(dz),\\
&=c^{-1}\int_{0}^{\infty}\mu(dv)\int^{0}_{-\infty} \mu(du)(e^{\beta v}-e^{\beta u})(-uv-\frac{\beta(u+v)}{2\sqrt{m}}),
\end{flalign*}
the last equality is followed from the construction of $\mu^{(m)}$ in (\ref{distri}).
Similarly, we can calculate that
\begin{flalign*}
b_{m}(x)
&=m\int_{|y-x|\leq1}(y-x)\Pi_{m}(x,dy),\\
&=m\int_{\frac{1}{\sqrt{m}}z\leq1}\frac{1}{\sqrt{m}}z\mu^{(m)}(dz),\\
&=-\frac{\beta}{2},
\end{flalign*}
and
\begin{flalign*}
\Delta_{m}^{\varepsilon}
&=m \Pi_{m}(x,\r \backslash B(x,\varepsilon))=\int_{x+\frac{1}{\sqrt{m}}z \in \r \backslash B(x,\varepsilon)}\mu^{(m)}(dz),\\
&=\int_{\frac{1}{\sqrt{m}}z \in \r \backslash B(0,\varepsilon)}\mu^{(m)}(dz)=\int_{|\frac{1}{\sqrt{m}}z|> \varepsilon}\mu^{(m)}(dz)=0.
\end{flalign*}
obviously, which satisfy the conditions (\ref{c1}), (\ref{c2}), (\ref{c3}).

\

(2){\it relatively compact of $\{H^{(m)}(\cdot)\}_{m\geq1}$}.

With  (\ref{generator}) in hand, by Corollary 4.3 of Chapter 4 in (\cite{Kz}), we need only to verify the following condition, that is, for all $\epsilon , T > 0$ there exist compact sets $K_{\epsilon, T} \subset \r$ such that
\begin{equation}\label{tight}
\inf_{m}P\{ H^{(m)}(t) \in K_{\epsilon, T} ~for~ all~ t\leq T\} \geq 1 - \epsilon.
\end{equation}
At first recall that for each $m\geq 1$, $X_{i}^{(m)}, i\geq 1$ are i.i.d. random variables, determined by  $\mu^{(m)}$ in (\ref{distri}). It is easy to calculate the moments of  $X_{i}^{(m)}$,
\begin{flalign*}
EX_{i}^{(m)} =&c^{-1}\int_{0}^{\infty}\mu(dv)\int^{0}_{-\infty}\mu(du)(e^{\beta v}-e^{\beta u}) \cdot \nonumber \\
&\left[ \left(\frac{-v}{u-v} - \frac{\beta}{2\sqrt{m}(u-v)} \right) u + \left(\frac{u}{u-v} + \frac{\beta}{2\sqrt{m}(u-v)} \right) v \right]\\
=& -\frac{\beta}{2\sqrt{m}},
\end{flalign*}
for each $i\geq 1$, and similarly
\begin{flalign}\label{var}
D(X_{i}^{(m)})= c^{-1}\int_{0}^{\infty}d F(v)\int^{0}_{-\infty}d F(u)(e^{\beta v}-e^{\beta u})(-uv-\frac{\beta(u+v)}{2\sqrt{m}}) -\frac{\beta^{2}}{4m}.
\end{flalign}

Let $\eta_{i}^{(m)}:=(X_{i}^{(m)}-EX_{i}^{(m)})/\sqrt{m}=(X_{i}^{(m)}+\frac{\beta}{2\sqrt{m}})/\sqrt{m}$ and $M_{n}^{(m)}=\sum_{i=1}^{n}\eta_{i}^{(m)} $, then for each $m\geq1$, $\eta_{i}^{(m)}, i\geq 1$ are i.i.d. random variables with $E\eta_{i}^{(m)}=0$, $Var \eta_{i}^{(m)} = \frac{DX_{i}^{(m)}}{m}$.

\begin{flalign}\label{sup}
P\left\{\sup_{0\leq x \leq T}|H^{(m)}(x)| \geq \lambda\right\}
&= P\left\{\sup_{1 \leq n \leq \lceil mT \rceil }\frac{1}{\sqrt{m}}|\sum_{i=1}^{n}X_{i}^{(m)}|\geq \lambda\right\},\nonumber\\
&= P\left\{\sup_{1 \leq n \leq \lceil mT \rceil }\left| \frac{1}{\sqrt{m}}\sum_{i=1}^{n}(X_{i}^{(m)}+\frac{\beta}{2\sqrt{m}}) - \frac{n\beta}{2m}\right| \geq \lambda\right\} \nonumber  \\
&\leq P\left\{\sup_{1\leq n \leq \lceil mT \rceil}|M_{n}^{(m)}|\geq \lambda \right\},
\end{flalign}
By Kolmogorov's maximal inequality (page 61, \cite{Dur}), (\ref{sup}) can be continuous
\begin{flalign*}
P\left\{\sup_{0\leq x \leq T}|H^{(m)}(x)| \geq \lambda\right\}
&\leq\frac{ Var M_{\lceil mT \rceil}^{(m)}}{\lambda^{2}}=\frac{\lceil mT \rceil \frac{DX_{i}^{(m)}}{m}}{\lambda^{2}}\leq \frac{(T+1)DX_{i}^{(m)}}{\lambda^{2}}\leq\frac{2TC}{\lambda^{2}},
\end{flalign*}
where $\lceil mT \rceil$ denotes the smallest integer greater than $mT$. The last inequality is because the sequence of
\{$DX_{i}^{(m)}$\}  is   convergent as $m\to\infty$ (see (\ref{var})), and then there is a constant $ C>0$ such that $DX_{i}^{(m)} \leq C$.  Now we choose $\lambda$ large enough such that $\frac{2TC}{\lambda^{2}}<\epsilon$, as a consequence (\ref{tight}) is proved if we  take $K_{\epsilon, T}=[-\lambda-1,\lambda+1]$.
\qed

\section{\it $\kappa$-transient random walk in random environment }           \label{s1.2:random}

Let us start from a $\kappa$-transient RWRE $\{Z_{n}\}_{n\geq0}$ with the i.i.d. environment series $\{\omega_i, i\in Z\}$ satisfying the condition (\ref{rwrec}), and assume the marginal distribution of $\omega_i$ is $\nu$. Recall $\rho_{i}=\frac{1-\omega_{i}}{\omega_{i}}$, $i \in \z$;  we write  the distribution of  $\log \rho_{0}$ as $\pi$. Now we can  construct a sequence of $\kappa$-transient RWRE as the following four steps,

\
$\clubsuit$ {\it step 1.   Define $\pi^{(m)}$. }

 We will construct a sequence of $\pi^{(m)}$ by (\ref{distri}) replacing $\mu=\pi$ and $\beta=\kappa$ respectively, i.e., for any continuous bounded function $\varphi$,
\begin{flalign}\label{1distri}
\int_{\r} \varphi(x) \pi^{(m)}(dx)  =& c^{-1}\int_{0}^{\infty}\pi(dv)\int^{0}_{-\infty}\pi(du)(e^{\kappa v}-e^{\kappa u}) \cdot \nonumber \\
&\left[ \left(\frac{-v}{u-v} - \frac{\kappa}{2\sqrt{m}(u-v)} \right) \varphi(u) + \left(\frac{u}{u-v} + \frac{\kappa}{2\sqrt{m}(u-v)} \right) \varphi(v) \right],
\end{flalign}
where $c:= \int_{0}^{\infty}(e^{\kappa x}-1)\pi(dx)=\int^{0}_{-\infty}-(e^{\kappa x}-1)\pi(dx) > 0$.

$\clubsuit$ {\it step 2.   Define $\{\delta_{i}^{(m)}\}_{i \in \z}$. }

 Let $\{\delta_{i}^{(m)}\}_{i \in \z}$ be a sequence of i.i.d. random variables with distribution $\pi^{(m)}$, which is defined in (\ref{1distri}).

$\clubsuit$ {\it step 3.   Define ``random environment" $\omega_{i}^{(m)}$. } Let
\begin{equation}\label{envir}
\omega_{i}^{(m)}:=\left( 1 + e^{\delta^{(m)}_{i}/\sqrt{m}} \right)^{-1}.
\end{equation}

$\clubsuit$ {\it step 4.   Define ``random walk in random environment" $Z^{(m)}$. }

 Denote $Z^{(m)}=\{Z_{n}^{(m)}\}_{n\geq0}$ as the random walk associated with the random environment $\omega^{(m)}=\{ \omega_{i}^{(m)}\}_{i \in \z}$, i.e.,
\begin{equation} \label{1rwre}
\p \left( Z_{n+1}^{(m)} = j ~|~Z^{(m)}_{n} = i; \omega^{(m)} \right) = \left\{\begin{array}{lll} \omega^{(m)}_{i}, & {\rm if}
\ j = i+1,
\\
 1-\omega^{(m)}_{i}, &  {\rm if} \ j = i-1,
\\
0, &{\rm otherwise}.
\end{array}
\right.
\end{equation}
From the above construction, we have

\begin{proposition}\label{p2} For each $m\geq1$,  the RWREs $Z^{(m)}=\{Z_{n}^{(m)}\}_{n\geq0}$ are $\kappa_m$-transient RWRE, i.e., each of them  satisfies condition (\ref{rwrec}) with some $\kappa_m>0$.
\end{proposition}\label{p2}
\noindent{\it Proof}~~
For each $m\geq1$, denote $\rho_{i}^{(m)}=\frac{1-\omega_{i}^{(m)}}{\omega_{i}^{(m)}}$, then, we have
\begin{equation*}
E\log \rho_{0}^{(m)}=-\frac{\kappa}{2m}<0,
\end{equation*}
so $Z^{(m)}$ is transient to $+\infty$. In addition, by the convex of the function of $f(s)=E\exp{s\log\rho_{0}^{(m)}},~s\geq0$ and $f'(0)=E\log\rho_{0}^{(m)}<0$, we know that $\exists~\kappa_{m}>0$, s.t. $E(\rho_{0}^{(m)})^{\kappa_{m}}=1$. So, for every $m\geq1$, RWRE $Z^{(m)}$ is an $\kappa_m$-transient RWRE, i.e., it  satisfies condition (\ref{rwrec}), a Kesten's type RWRE.\qed

\begin{theorem}\label{th2}
{\it Under the assumptions (\ref{rwrec}) and uniformly elliptic, i.e. $\exists ~0 < \varepsilon < 1/2$ such that $P(\omega_{0} \in [\varepsilon,1-\varepsilon])=1$; $Z^{(m)}$ be defined in (\ref{1rwre}). Then, as $m\to\infty$
\begin{equation}\label{RB}
\left\{ \frac{1}{m}Z^{(m)}_{[m^2 t]}, t\geq 0 \right\}\longrightarrow \{ X_{V}(t), t\geq 0\}
\end{equation}
in distribution in $\mathcal{D}[0,\infty)$. Where $X_{V}=(X_{V}(t), t\geq0)$ a diffusion process in random potential $ {V}(x)=\sigma \cdot W(x)-\frac{\kappa}{2}x$  (defined in (\ref{CX}) and (\ref{PV}) respectively ), and $\sigma^{2}$ defined in (\ref{sigma}) replacing $\mu=\pi$ and $\beta=\kappa$.
}
\end{theorem}

\noindent\begin{remark} Recall Seignourel \cite{S00} proved that the scaling limit of Sinai's random walk random environment (i.e. $E\log \rho_{0}=0$) is Brox's diffusion(i.e V=Wiener process).\end{remark}
\noindent{\bf  Theorem A} (Seignourel, \cite{S00}, 2000) { Assume $E\log \rho_{0}=0$, and denote $\sigma^{2}:= E( \log \rho_{0} )^{2}$. For every $m\geq1$ and for every $i \in \z$, let
\begin{equation}\label{envir1}
\omega_{i}^{(m)}:=( 1+(  \rho_{i} )^{\frac{1}{\sqrt{m}}})^{-1}=\left( 1 + e^{\log \rho_{i}/\sqrt{m}} \right)^{-1}.
\end{equation}
Then, as $m\to\infty$
\begin{equation*}
\left\{ \frac{1}{m}Z^{(m)}_{[m^2t]}, t\geq 0 \right\}\longrightarrow \{ X_{V}(t), t\geq 0\}
\end{equation*}
in distribution in $\mathcal{D}[0,\infty)$, where $\{ X_{V}(t), t\geq 0\}$ is the diffusion process with random potential $V=\sigma$ $\cdot$ Brownian motion.
}

{\it Comparing the ``environment" definition in (\ref{envir}) and (\ref{envir1}), our ``measure change" related to $m\geq 1$ is through the $\delta^{(m)}_{i}$ in (\ref{envir}). \qed }

\

%%%%%%%%%%%%%%%%%%%%%%%%%%%Proof Theorem1.2%%%%%%%%%%%%%%%%%%%%%%%

\noindent
{\it Proof of Theorem \ref{th2}} We can follows the proof of the Theorem 1 in \cite{S00} (Seignourel, 2000) almost line by line, except that the $W^m(x)$ in Theorem 12 of \cite{S00} is now that
\begin{equation}\label{rwrepoten}
V^{(m)}(x)= \left\{\begin{array}{lll} \frac{1}{\sqrt{m}}\sum_{i=1}^{[ mx ]} \delta_{i}^{(m)}, & {\rm if}
\ [ mx ]\geq 1,
\\
 0, &  {\rm if}\  0 \leq [ mx ]\leq 1,
\\
-\frac{1}{\sqrt{m}}\sum_{i=1}^{[ mx ]}    \delta_{i}^{(m)}, & {\rm if}
\ [ mx ]< 0,
\end{array}
\right.
\end{equation}
where $\{\delta_{i}^{(m)}\}_{i \in \z}$ be a sequence of i.i.d. random variables with distribution $\pi^{(m)}$, which is defined in (\ref{1distri}). By the definition of $\omega_{i}^{(m)}$ in (\ref{envir}), actually we know with $\rho_{i}^{(m)}=\frac{1-\omega_{i}^{(m)}}{\omega_{i}^{(m)}}$,
 $$\delta_{i}^{(m)}={\sqrt{m}}\log \rho_{i}^{(m)}.$$
The proof will complete if we can verify the sequence of $V^{(m)}(x)$ convergence weakly to $ {V}(x)=\sigma \cdot W(x)-\frac{\kappa}{2}x$  (defined in  (\ref{PV}) ).

To this end, we need only to take $ X_{i}^{(m)}:= \delta_{i}^{(m)}$ in  Proposition \ref{th1}; and by the construction of the sequences of $\{ X_{i}^{(m)}\}$ and $\{\delta_{i}^{(m)}\}$,  we know that  $\{ X_{i}\}:=\{\log \rho_{i}\}$. The conditions (\ref{rwrec}) and uniformly elliptic conditions enable the sequences of $ X_{i}^{(m)}$ (i.e., $\{\delta_{i}^{(m)}\}$) fulfill the  conditions (\ref{rwc}), as a consequence the sequence of $V^{(m)}(x)$ convergence weakly to $ {V}(x)=\sigma \cdot W(x)-\frac{\kappa}{2}x$  by Proposition \ref{th1}. \qed


\begin{thebibliography}{99}
\def\nobibitem#1\par{}
\bibitem{BN12}
Borisov, I.~S. and Nikitina, N.~N. \ \
The distribution of the number of crossings of a strip by paths of the simplest random walks and of a wiener process with drift.
{\it Theory Probab. Appl.}, {\bf 56}(2012), pp. 126-132.

\bibitem{B86}
Brox, T.\ \
A one-dimensional diffusion process in a Wiener medium.
 {\it Ann.\ Probab.} {\bf 14 } (1986), 1206-1218.

\bibitem{C00}
Comets, F., Gantert, N. and Zeitouni, O. \ \
Quenched, annealed and functional large deviations for one-dimensional random walk in random environment.
{\it Probab. Th. Ralated Fields} {\bf 118} (2000), 65-114.




\bibitem{DPZ96}
  Dembo,A., Peres,Y. and  Zeitouni,O., Tail estimates for one-dimensional random walk in random environment,
{\it Communications in Mathematical Physics,} {\bf 181} (1996), 667-683.

\bibitem{Dur}
Durrett, R.\ \
 {Probability: Theory and Examples, 3rd Edition.}
Duxbury, 2004.

\bibitem{ESZ09}
 Enriquez,N.,  Sabot,C. and Zindy, O., Limit laws for transient random walks in random environment on $\ mathbb {Z} $,
{\it Annales de l'institut Fourier,} {\bf 59} (2009), 2469-2508.




\bibitem{Kurtz}
Ethier, S.~N.\ and Kurtz, T.~G.\ \
 {Markov Processes: Characterization and Convergence, Second edition.}
Wiley Series in Probability and Statistics, 2005.




\bibitem{G94}
Greven, A. \ and den Hollander, F. \ \
Large deviations for a random walk in random environment.
{\it Ann. Probab.} {\bf 22} (1994), 1381-1428.






\bibitem{H99}
Hu, Y., Shi, Z. and Yor, M. \ \
Rates of convergence of diffusions with drifted Brownian potentials.
{\it Trans. Amer. Math. Soc.} {\bf 351} (1999), 3915-3934.


\bibitem{K97}
Kawazu, K. and Tanaka, H. \ \
A diffusion process in a Brownian environment with drift.
{\it J. Math. Soc. Japan} {\bf 49} (1997), 189-211.

%\bibitem{K98}
%Kawazu, K. and Tanaka, H. \ \
%Invariance principle for a Brownian motion with large drift in a white noise environment.
%{\it Hiroshima Math. J.} {\bf 28} (1998), 129-137.

\bibitem{K75}
Kesten, H., Kozlov, M.V. and Spitzer, F. \ \
A limit law for random walk in random environment.
{\it Compositio Math.} {\bf 30} (1975), 145-168.

\bibitem{Kz}
Kurtz, T.~G.\ \
Approximation of Population Processes,
Society for industrial and applied mathematics, 1981.



\bibitem{S85}
Schumacher, S. \ \
Diffusions with random coefficients.
{\it Contemp. Math.} {\bf 41} (1985), 351-356.


\bibitem{S00}
Seignourel, P.\ \
Discrete schemes for processes in random media.
 {\it Probab.\ Theory Relat.\ Fields.} {\bf 118 } (2000), 293-322.


\bibitem{V09}
Stroock, D. W., Varadhan, S.R.S., \ \
{Multidimensional diffusion processes,}
{\it Springer} 2009


\bibitem{T01}
Taleb, M. \ \
Large deviations for a Brownian motion in a drifted Brownian potential.
{\it Ann. Probab.} {\bf 29} (2001), 1173-1204.


\bibitem{T95}
Tanaka, H. \ \
Diffusion processes in random environments.
{\it Proc. ICM(S.D. Chatterji, ed.)} 1047-1054. Birkh$\ddot{a}$user, Basel.




\bibitem{ze04}  Zeitouni, O. (2004). Random walks in random environment. \textit{LNM 1837, J. Picard (Ed.), pp. 189-312, Springer-Verlag Berlin Heidelberg}.

\end{thebibliography}
\end {document}